\documentclass[12pt]{amsart}
\usepackage{amsmath,amsthm,amsfonts,amssymb,eucal}
\usepackage[dvips]{graphics}

\newcommand{\Om}{\Omega}

\newcommand{\oq}{\ {\raise 7pt\hbox{${\scriptstyle\circ}$}}
\kern -7pt{
\hbox{$Q$}}}

\newcommand{\R}{ \mathbb R}
\newcommand{\N}{ \mathbb N}
\newcommand{\al}{ \alpha}

\newcommand {\bs}{\begin{split}}

\newcommand {\De}{\Delta}
\newcommand {\di}{\partial}
\newcommand {\ps}{\pi^2}
\newcommand {\psf}{\frac{\pi^2}{4}}

\newcommand{\1}
{{\,\vrule depth3pt height9pt}{\vrule depth3pt height9pt}
{\vrule depth3pt height9pt}{\vrule depth3pt height9pt}\,}

\newtheorem{thm}{Theorem}[section]
\newtheorem{cor}[thm]{Corollary}
\newtheorem{lem}[thm]{Lemma}

\theoremstyle{definition}

\theoremstyle{remark}
\newtheorem{rem}[thm]{Remark}

\numberwithin{equation}{section}
\newcommand{\be}{\begin{equation}}
\newcommand{\bd}{\begin{displaymath}}
\newcommand{\ed}{\end{displaymath}}

\newcommand{\bee}{\begin{equation}}
\newcommand{\ene}{\end{equation}}
\newcommand{\bes}{\begin{split}}
\newcommand{\ens}{\end{split}}

\newcommand{\bet}{\begin{thm}}
\newcommand{\ent}{\end{thm}}
\newcommand{\bel}{\begin{lem}}
\newcommand{\enl}{\end{lem}}
\newcommand{\bec}{\begin{cor}}
\newcommand{\enc}{\end{cor}}
\newcommand{\bep}{\begin{proof}}
\newcommand{\enp}{\end{proof}}
\newcommand{\ber}{\begin{rem}}
\newcommand{\enr}{\end{rem}}

\newcommand{\la}{\lambda}
\newcommand{\de}{\delta}
\newcommand{\Z}{\mathbb Z}

\makeatletter
\def\square{\RIfM@\bgroup\else$\bgroup\aftergroup$\fi
  \vcenter{\hrule\hbox{\vrule\@height.6em\kern.6em\vrule}\hrule}\egroup}
\makeatother

\begin{document}
\title[Thick obstacle
(\the\day.\the\month.\the\year)] {Trapped modes in a waveguide
with a thick obstacle}
\author[H. Hawkins]
{Helen Hawkins}
\address{Centre for Mathematical Analysis and Its Applications\\
University of Sussex\\ Falmer,  Brighton\\ BN1 9QH, UK}
\email{h.l.hawkins@sussex.ac.uk}
\author[L.Parnovski]{Leonid Parnovski}
\address{Department of Mathematics\\
University College London\\ Gower Street\\ London WC1E 6BT, UK}
\email{Leonid@math.ucl.ac.uk}
\date{\today}

\maketitle

\section{Introduction} 
\begin{figure}[ht]
\begin{center}
\scalebox{0.6}{\includegraphics{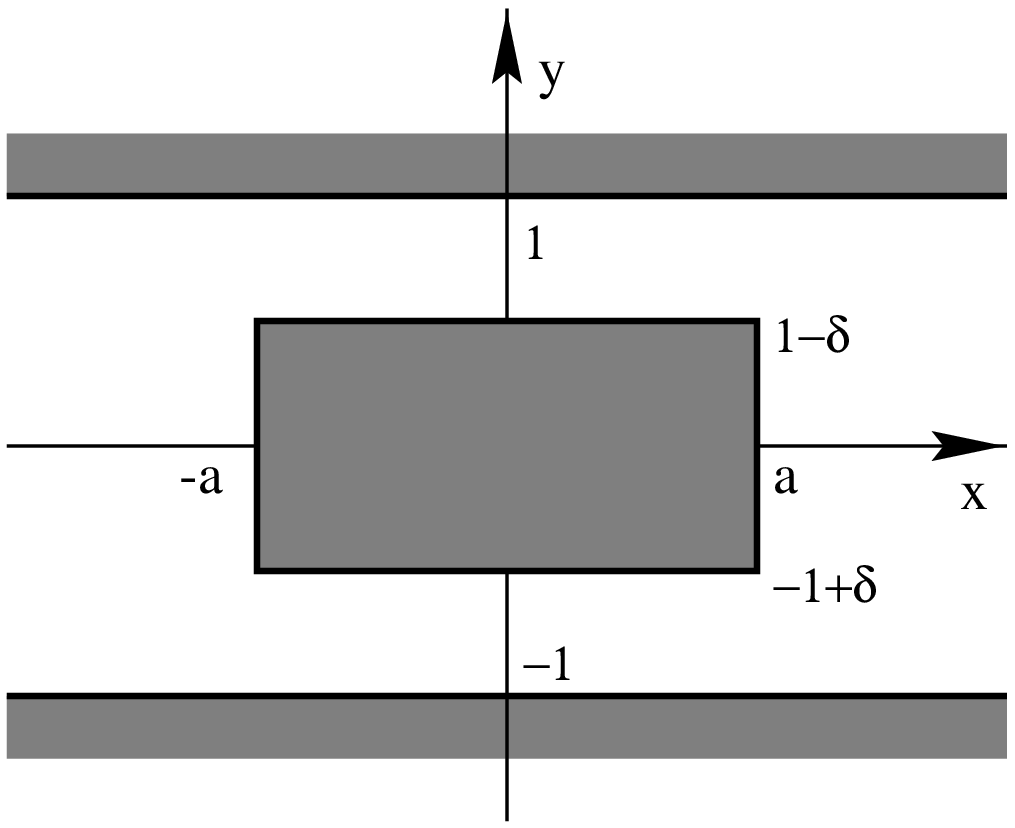}}
\end{center}
\caption{}
\label{fig_rectangle}
\end{figure}

The problem of finding necessary and
sufficient conditions for the existence of trapped modes in
waveguides has been known since 1943, \cite{Rel}. The problem is
the following: consider an infinite strip $M$ in $\R^2$ (or an
infinite cylinder with the smooth boundary in $\R^n$). The
spectrum of the (positive) Laplacian (with either Dirichlet or
Neumann boundary conditions) acting on this strip is easily
computable via the separation of variables; the spectrum is
absolutely continuous and equals $[\nu_0,+\infty)$. Here, $\nu_0$
is the first {\em threshold}, i.e. eigenvalue of the cross-section
of the cylinder (so $\nu_0=0$ in the case of Neumann conditions).
Let us now consider the domain $\tilde\Om$ (the {\em waveguide})
which is a smooth compact perturbation of $M$ (for example, we
insert an obstacle in $M$). The essential spectrum of the
Laplacian acting on $\tilde{\Om}$ still equals $[\nu_0,+\infty)$,
but there may be additional eigenvalues (so-called {\em
trapped modes}; the number of these trapped modes can be quite
large, see examples in \cite{Wit} and \cite{Par}). So, the problem
is in finding conditions for the existence or absence of such
eigenvalues and studying them when they exist. It is customary to
distinguish between two situations: the Dirichlet boundary
conditions (corresponding to the so-called {\em quantum
waveguides}) and the Neumann boundary conditions (corresponding to
the {\em acoustic waveguides}). In the Dirichlet case the first
threshold $\nu_0>0$, so the eigenvalues can occur outside
the essential spectrum. Such eigenvalues (not embedded into the
essential spectrum) are stable under small perturbations, and thus
they occur in a wide range of situations (see e.g. \cite{Ex}). On
the contrary, in the Neumann case any eigenvalue is
embedded into the continuous spectrum and is very unstable.
Therefore, the existence of such eigenvalues is usually due to
some symmetry (obvious or hidden) of the situation. In
\cite{EvaLevVas} it was shown that if the obstacle is symmetric
about the axis of the strip, then for a wide range of obstacles
there is (at least one) eigenvalue. Later, in \cite{DavPar} more
conditions, necessary as well as sufficient for the existence of
eigenvalues were established. Also, in that paper the example of
a hidden symmetry resulting in the existence of an eigenvalue
was given. Once the existence of eigenvalues is established, it is
natural to ask how many of them there are and how they behave. In
the case of a symmetric obstacle the problem splits into two
problems on the half of $\tilde\Om$ (obtained by cutting
the initial waveguide
along its axis); one of these problems corresponds to the
additional Dirichlet conditions on the added boundary; the other
has additional Neumann conditions. The essential spectrum of the
first sub-problem starts at the first non-zero threshold $\nu_1>0$
(in the case of the strip of width $2$ which we will consider in
this paper, $\nu_1=\pi^2/4$), and the second sub-problem still has
essential spectrum growing from zero (see next section for more
details). Paper \cite{KhaParVas} studied what happens if the
symmetric obstacle becomes long (in the direction of the axis of
the strip). It is proved there that the number of eigenvalues of
the Dirichlet sub-problem below $\nu_1$ is of the order of the
length of the obstacle. In the present paper we study another
regime of the asymptotic behaviour of such eigenvalues: suppose
that the obstacle is a rectangle placed symmetrically on the axis
of the strip (see figure \ref{fig_rectangle}).
Let the width of the strip be $2$, the
length of the rectangle (in the direction of the axis of the
strip) be $2a$ and the distance from the rectangle to the sides of
the waveguide (in the direction
orthogonal to the axis) be $\de$. When $\de=0$, the domain
degenerates to the union of two semistrips. We are
interested in the behaviour of the eigenvalues which lie below the
first non-zero threshold $\nu_1=\pi^2/4$ when $\de\downarrow 0$,
in particular, the rate at which they tend to the threshold.

\begin{figure}[ht]
\begin{center}
\scalebox{0.6}{\includegraphics{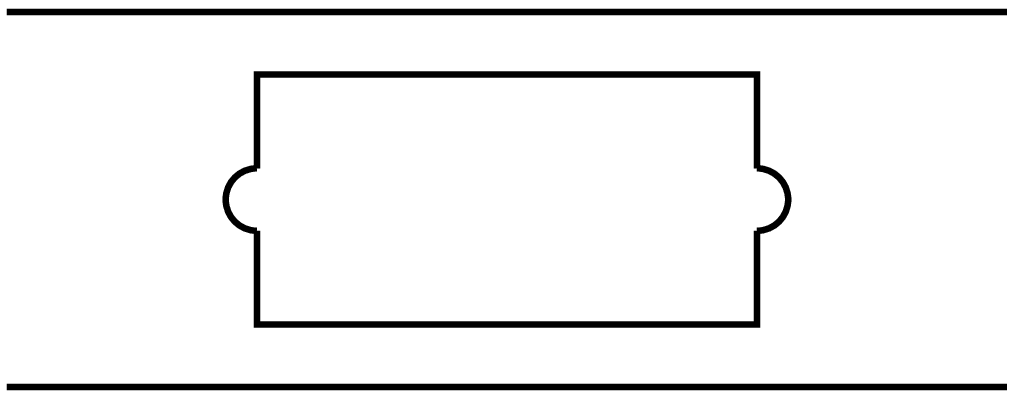}}
\end{center}
\caption{}
\label{fig_bump}
\end{figure}

\begin{figure}[ht]
\begin{center}
\scalebox{0.6}{\includegraphics{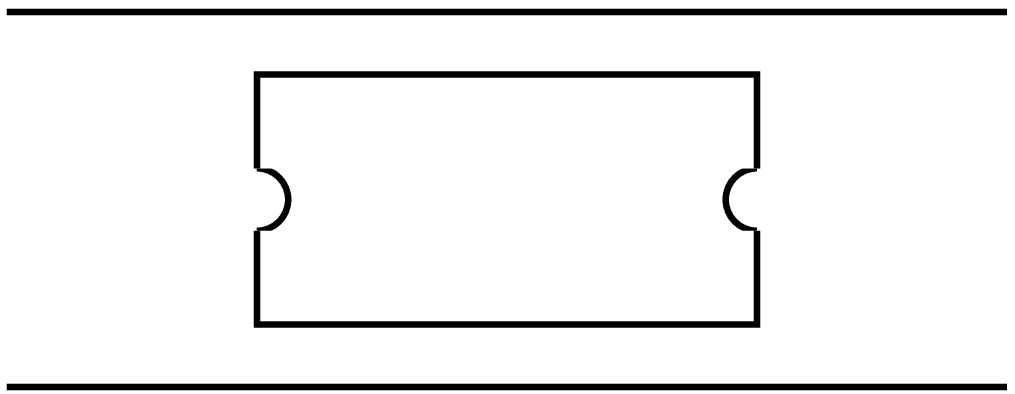}}
\end{center}
\caption{}
\label{fig_indent}
\end{figure}

The choice of the rectangle as an obstacle is motivated by the
following considerations: suppose for simplicity that $0<a<1$. Then,
if the obstacle has the shape as in figure \ref{fig_bump},
the (unique) eigenvalue
will stay away from the threshold $\nu_1$ (this can be proved using the
same method as in \cite{EvaLevVas}). On the other hand, if the
obstacle has the shape as in figure \ref{fig_indent},
it was shown in \cite{DavPar} that
for small enough $\de$ there are no eigenvalues below $\nu_1$ at all. The
case of a rectangle is an intermediate one:
for any $\de$ there is a unique eigenvalue
which converges towards $\nu_1$ as $\de\to 0$. This makes
the case of a rectangular obstacle such an interesting one. A slightly
different problem about the rate of convergence of an eigenvalue to a
threshold (in the context of a quantum waveguide) was considered
in \cite{ExnVuh} and \cite{Pop}, and that problem turned out to be
quite difficult (so that one has to work a lot even to get the correct
order of convergence). Our problem, on the contrary,
is relatively easy, and one can get the first asymptotic term
without too much difficulty (in fact, only the first transversal
mode contributes towards the first asymptotic term). We think
that one can also obtain the second asymptotic term (by studying
further transversal modes), but we have not done this in our
paper. The result we have obtained is rather surprising in the
sense that the rate of convergence of an eigenvalue towards the
threshold depends on whether $a$ is an integer or not. We postpone
the precise formulation of the result until the next section.

The
proof of our result will go along the standard lines. To
estimate the eigenvalue from above, we will produce the
test-function (or functions, if we have several eigenvalues). To
obtain the precise asymptotic constant, the test-function has to
be chosen with great care. In order to estimate the eigenvalue
from below, we use the technique of estimating quadratic forms,
similar to the method of transference of excess energy (see
\cite{DavPar}). There is a small difference between our
approach and the method of \cite{DavPar}. Namely, instead of
comparing the integrals of the function along different sub-regions
of $\tilde{\Om}$ (which was the key tool in \cite{DavPar}),
in our paper we compare such integrals with values of the function
in certain points.

The rest of the paper is organised in the following way: in the
next section we give some preliminary information and formulate
the main theorem \ref{main}; sections
3-5 are devoted to the proof of this theorem. For the
convenience of the reader we discuss first (in the section 3)
the easiest case when $a<1$. Section 4 deals with the case $a=1$
(so that in both these sections we have only one eigenvalue).
Finally, in section 5 we explain which changes should be made in
the proof when $a$ is arbitrary (and there are several
eigenvalues).

{\bf Acknowledgement} We are very grateful to D. Vassiliev and M. Levitin
for helpful suggestions and comments.

\section{Preliminaries}
We consider the domain $\tilde\Om=(-\infty,\infty)\times (-1,1)\setminus
(-a,a)\times (-(1-\de),1-\de)$, $0<a$, $0<\de<1$. The spectrum of $-\De$ with
Neumann boundary conditions on $\tilde\Om$ is the interval
$[0,+\infty)$. To make the study of the eigenvalues easier, we
split $L^2(\tilde{\Om})$ into several subspaces invariant with respect of
the action of $-\De$. First, let $\Om$ be the half of $\tilde\Om$:
$\Om=(-\infty,\infty)\times (0,1)\setminus
(-a,a)\times (0,1-\de)$. It is well known (see \cite{EvaLevVas}) that
if we consider the operator $L$ which acts as $-\De$ on $\Om$ with
Dirichlet boundary conditions on $\{y=0\}$ and Neumann boundary
conditions on the rest of the boundary $\di\Om$, then eigenvalues
of $L$ are at the same time eigenvalues of the Neumann Laplacian
on $\tilde\Om$. Moreover, since the essential spectrum of $L$ is
$[\pi^2/4,+\infty)$, we can study eigenvalues below
$\nu_1:=\pi^2/4$ using the variational approach. It is also
convenient to make the further reduction of the domain and
consider two problems on $\Om':=(0,\infty)\times (0,1)\setminus
(0,a)\times (0,1-\de)$: one problem, called $L_D$, has Dirichlet
conditions on $\{x=0\}\cup\{y=0\}$ and Neumann conditions
elsewhere; the other problem, called $L_N$, has Dirichlet
conditions on $\{y=0\}$ and Neumann conditions
elsewhere (see \cite{KhaParVas} for more details of this
decomposition). Then $\sigma(L)=\sigma(L_D)\cup\sigma(L_N)$, i.e.
the spectrum of $L$ is the union of spectra of $L_D$ and $L_N$.
Let $\la_1\le\dots\le\la_n$ be the eigenvalues of
$L$ lying below $\nu_1$. Using the approach of \cite{KhaParVas},
together with the test-function from \cite{EvaLevVas}, it is easy
to show that $n=-[-a]$ (this will also follow from the proof of our
main theorem). Moreover, if $n$ is even, then half of these
eigenvalues come from $L_N$, and another half comes from $L_D$. If
$n$ is odd, then the spare eigenvalue is due to $L_N$. It is also
known that the eigenvalues coming from $L_N$ and $L_D$ are alternating
and that $\la_1\in\sigma(L_N)$.
Thus, the top eigenvalue $\la_n$ is an eigenvalue of $L_N$ if and only if $n$
is odd. If we fix $a$ and let $\de\to 0$, then all but the last
eigenvalue remain bounded away from $\nu_1$, i.e.
$\la_{n-1}<C<\nu_1$ uniformly over $\de$ (see
\cite{KhaParVas}). On the other hand, $\la_n\to\nu_1$ as $\de\to
0$. 
Now we can formulate our main result.

\begin{thm}\label{main}
If $a\not\in\N$, then
\bee\label{main1}
\la_n=\psf-M(a)\de^2+O(\de^3)
\end{equation}
as $\de\to 0$, where
\bee\label{M1}
M(a)=\ps\tan^2\left(\frac{\pi\{a\}}{2}\right)= \left\{ \begin{array}{ll}
\ps \tan^{2} \left(\frac{\pi a}{2}\right) & \textrm{$ 2l<a<2l+1,$}\\
\ps \cot^{2} \left(\frac{\pi a}{2}\right) & \textrm{$2l+1<a<2l+2$.}
\end{array} \right.
\end{equation}
Here $\{a\}:=a-[a]$ is the fractional part of $a$.
If $a\in\N$, then
\bee\label{main2}
\la_n=\psf-M(a)\de^{2/3}+O(\de^{4/3})
\end{equation}
as $\de\to 0$, where
\bee\label{M2}
M(a) = \left(\frac{\ps}{a}\right)^{2/3}.
\end{equation}
\end{thm}

The rest of the paper is devoted to the proof of this theorem.

\section{$a<1$}\label{a<1} In this case there is only one eigenvalue (coming
from $L_N$), and we denote this eigenvalue by $\la$. Also, $M(a)=
\pi^2\tan^{2}\bigl(\frac{\pi a}{2}\bigr)$.
Obviously, \eqref{main1} is equivalent to the
following two inequalities:
\bee\label{main1a}
\la\le \frac{\pi^2}{4}-M(a)\de^2+O(\de^3)
\end{equation}
and
\bee\label{main1b}
\la\ge \frac{\pi^2}{4}-M(a)\de^2-O(\de^3).
\end{equation}
The strategy of the proof will be quite standard for problems
of this sort: to prove \eqref{main1a}, we will construct the
test-function $\phi$ satisfying Dirichlet conditions at $\{x=0\}$
for which the Rayleigh quotient
\bee\label{Rq}
Q(\phi):=\frac{\int_{\Om'}|\nabla\phi|^2 dxdy}{\int_{\Om'}|\phi|^2 dxdy}=
\frac{\pi^2}{4}-M(a)\de^2+O(\de^3),
\end{equation}
and to prove \eqref{main1b}, we will estimate the
quadratic form of $L_{N}$ from below. It is relatively easy to
construct the test-function $\phi$ for which
\bee\label{test1}
Q(\phi)=\frac{\pi^2}{4}-K\de^2-O(\de^3),
\end{equation}
but the constant $K$ is worse than $M(a)$. For example, let us
denote by
\bee\label{eta}
\eta_k(y):=\sqrt 2 \sin\frac{(2k-1)\pi y}{2}
\end{equation}
($k\ge 1$)
the normalized eigenfunctions of the cross-section of the strip at
infinity. Then if we choose
\bee\label{phi1}
\phi(x,y)=\begin{cases}
\eta_1(y)e^{-\de (x-a)}, & x\ge a\\
\eta_1(y), & x\le a,
\end{cases}
\end{equation}
this function will satisfy \eqref{test1} with $K=\pi^2a-1<M(a)$. In
order to get the precise constant, we have to correct the function
\eqref{phi1} in the region above the obstacle. This correction is
not obvious, and in order to understand it, we will first prove
\eqref{main1b}. To begin with, we decompose $\Om'$ into two parts:
\bee\label{Oml}
\Om_l:=(0,a)\times (1-\de,1)
\end{equation}
and
\bee\label{Omr}
\Om_r:=(a,+\infty)\times (0,1).
\end{equation}
The estimate of the quadratic form will be
different in $\Om_l$ and $\Om_r$. In each of these regions we will
use certain one-dimensional results to obtain the estimates. The
first lemma is rather trivial; it will take care of $\Om_r$.

\bel\label{right}
For any function $\varphi \in C^{1}(\mathbf{R^{+}})$ and any $m>0$
we have:
\begin{equation}\label{resultr}
 m \varphi^{2}(0) \leq \int_{0}^{\infty}
 \left( \varphi'(x)^{2} + m^{2} \varphi^{2}(x) \right) dx.
\end{equation}
with equality only when $\varphi(x)=Ce^{-m x}$.
\enl

\bep
Taking into account that
\begin{equation*}
\varphi^{2}(0)=-\int_{0}^{\infty} \frac{ d}{dx} (\varphi^{2}(x))
dx=-2\int_{0}^{\infty}  \varphi(x) \varphi'(x) dx,
\end{equation*}
we see that \eqref{resultr} is equivalent to
\begin{equation*}
\int_{0}^{\infty} (\varphi'(x)+ m \varphi(x))^2 dx\ge 0,
\end{equation*}
which makes both statements of lemma obvious.
\enp
The second lemma will help us to deal with $\Om_l$; this lemma is
slightly more subtle.

\begin{lem}\label{left}
Let $0<a<1$ and $f \in C^{2}[0,a], \, f'(0)=0$. Then
\be\label{leftleft}
\int_{0}^{a} \left( \frac{df}{dx}\right)^{2}dx -
\psf \int_{0}^{a} f^{2}(x) dx + \frac{\sqrt{M(a)}\, f^{2}(a)}{2}  \geq 0,
\end{equation}
where
\be\label{aless1_A}
M(a) = \ps \tan^{2} \left( \frac{\pi a}{2}\right).
\end{equation}
Moreover, equality is reached in \eqref{leftleft} if
\be\label{aless1_equality}
f(x) = C \cos\left(\frac{\pi x}{2}\right).
\end{equation}
\end{lem}
\begin{proof}
Without loss of generality we can assume $f$ satisfies the extra boundary condition
\be\label{bc}
f'(a) = - \frac{\sqrt{M(a)}\, f(a)}{2}.
\end{equation}
Indeed, given any function $f$ we can choose another function $g$
such that $g$ satisfies \eqref{bc} and the difference of the left
hand sides of \eqref{leftleft} for $f$ and $g$ is arbitrarily small.

We now note that
\be\label{aless1_quadratic_form}
\int_{0}^{a} \left(\frac{df}{dx}\right)^{2} dx + \frac{\sqrt{M(a)}\,f^{2}(a)}{2}
\end{equation}
is the quadratic form of $-\frac{d^{2}}{dx^{2}}$ with the boundary conditions
$f'(0)=0, f'(a)=-\frac{\sqrt{M(a)}\,f(a)}{2}$.

The eigenvalues of this operator are the values of $\mu^{2}$,
where $\mu$ satisfy
\be\label{aless1_M(a)_tan}
\mu \tan(\mu a) = \frac{\sqrt{M(a)}}{2},
\end{equation}
corresponding eigenfunctions being $\cos(\mu x)$. Therefore, the first
eigenvalue equals $\psf$ if and only if $M(a)$ is given by
\eqref{aless1_A}. This finishes the proof of the lemma.
\end{proof}

Now we will prove \eqref{main1b}. To do this, it is enough to show
that whenever $u\in C^{\infty}(\Om')$, $u(x,0)=0$, $\frac{\di u}{\di
x}(0,y)=0$, $\frac{\di u}{\di y}(x,1)=0$
the following inequality is satisfied:
\bee\label{quadratic}
\int\!\!\!\int_{\Omega'} \vert \nabla u |^{2} dxdy - \left(
\frac{\pi^{2}}{4} - M(a)\delta^{2}-O(\de^3) \right) \int\!\!\!\int_{\Omega'}
|u|^{2} dxdy\ge 0
\end{equation}
where $M(a)=\ps\tan^2 \left(\frac{\pi a}{2}\right)$. Let us decompose
the LHS of
\eqref{quadratic} as
\bee
\bes
&\left(\int\!\!\!\int_{\Omega_l} \vert \nabla u |^{2} dxdy -
\frac{\pi^{2}}{4}   \int\!\!\!\int_{\Omega_l}
|u|^{2} dxdy\right)\\
+&\left(\int\!\!\!\int_{\Omega_r} \vert \nabla u |^{2} dxdy -
\frac{\pi^{2}}{4} \int\!\!\!\int_{\Omega_r}
|u|^{2} dxdy\right)\\
+ &\left(M(a)\delta^{2}+O(\de^4) \right)\int_{\Omega}
|u|^{2} dxdy\\
=&:L(u)+R(u)+ \left(M(a)\delta^{2}+O(\de^4) \right)\int_{\Omega}
|u|^{2} dxdy
\end{split}
\end{equation}
(the definitions of $\Om_l$ and $\Om_r$ are given in \eqref{Oml}
and \eqref{Omr}).
The main idea of the proof is the following: it is obvious that
$R(u)>0$; moreover, there is a certain extra amount of energy in
$\Om_r$ to spare. We wish to transfer this excess of energy into
$\Om_l$ using the information of one-dimensional problems,
similar to the approach in \cite{DavPar}. However, if we
do it precisely in the same way as in \cite{DavPar}, the estimate
we obtain will be too rough. Therefore, instead we transform the
excess of energy of $u$ over $\Om_r$ into an extra positive term
involving the values of $u$ on the boundary between $\Om_r$ and
$\Om_l$. To be more precise, we will show that
\bee\label{Omr1}
R(u)\ge \frac{\sqrt{M(a)}}{2}\int_{1-\de}^1 |u(a,y)|^2 dy
\end{equation}
and
\bee\label{Omr2}
L(u)+\frac{\sqrt{M(a)}}{2}\int_{1-\de}^1 |u(a,y)|^2 dy\ge 0,
\end{equation}
which obviously would lead to \eqref{quadratic}.

We start by examining $R(u)$ in more details. Since $u(x,0)=0$
and $\frac{\di u}{\di y}(x,1)=0$,
we can decompose $u$ in the
Fourier series when $x\ge a$:
\bee\label{Fourier2}
u(x,y) = \sum_{k=1}^{\infty} u_{k} (x) \eta_{k}(y)
\end{equation}
($\eta_j$ is defined in \eqref{eta})
After simple computations we obtain:
\begin{equation*}
\int\!\!\!\int_{\Omega_{r}} \vert \nabla u \vert^{2}=
\sum_{k=1}^{\infty} \int_{a}^{\infty} \left( u_{k}'(x)
\right)^{2} dx + \frac{\pi^{2}}{4} \sum_{k=1}^{\infty} (2k-1)^{2}
\int_{a}^{\infty}u_{k}^{2}(x) dx
\end{equation*}
and
\begin{equation*}
\int\!\!\!\int_{\Omega_{r}} |u|^{2} dxdy = \sum_{k=1}^{\infty}
\int_{a}^{\infty} u_{k}^{2}(x) dx.
\end{equation*}
These formulae imply
\begin{equation}\label{Omr31}
\bes
R(u)&=\sum_{k=1}^{\infty} \int_{a}^{\infty} \left(u_{k}'(x)
\right)^{2} dx + \frac{\pi^{2}}{4} \sum_{k=1}^{\infty} (2k-1)^{2}
\int_{a}^{\infty} u_{k}^{2}(x) dx\\
& - \frac{\pi^{2}}{4}
\sum_{k=1}^{\infty} \int_{a}^{\infty} u_{k}^{2}(x) dx
+ M(a) \delta^{2} \sum_{k=1}^{\infty} \int_{a}^{\infty} u_{k}^{2}
(x)\\
&\geq  \sum_{k=1}^{\infty} \int_{a}^{\infty} \left(u_{k}'(x)
\right)^{2} dx + M(a)\delta^{2}\sum_{k=1}^{\infty}  \int_{a}^{\infty}
u_{k}^{2}(x) dx
\end{split}
\end{equation}
when $\de$ is small enough. Now lemma \ref{right} implies
\begin{equation}\label{Omr3}
R(u)\ge  \sqrt{M(a)} \, \delta \sum_{k=0}^{\infty}
u_{k}^{2} (a).
\end{equation}
On the other hand, the RHS of \eqref{Omr1} is
\bee\label{Omr4}
\bes
&\frac{\sqrt{M(a)}}{2}\int_{1-\de}^1 |u(a,y)|^2 dy\\&=
\frac{\sqrt{M(a)}}{2}
\sum_{k=1}^{\infty}u_{k}^{2} (a) \int_{1-\delta}^{1} 2 \sin^{2}
\left( \frac{\pi}{2} (2k-1) y \right) dy\\
&=\frac{\sqrt{M(a)}}{2}
\sum_{k=1}^{\infty}u_{k}^{2} (a)
\bigl(\delta + \frac{\sin \left(\pi \delta (2k-1) \right)}{\pi
(2k-1)}\bigr)\le \sqrt{M(a)}\de \sum_{k=1}^{\infty}u_{k}^{2} (a),
\end{split}
\end{equation}
which, together with \eqref{Omr3}, proves \eqref{Omr1}. Equation
\eqref{Omr2} follows immediately if we apply lemma \ref{left} to
the function $u(\cdot,y)$ and then integrate the result over $y$.
This finishes the proof of \eqref{main1b}.

In order to construct
the test-function satisfying \eqref{Rq}, we try to change all
inequalities in the proof of the lower bound into equalities. In
other words, we need to have equalities in \eqref{Omr2}, \eqref{Omr3},
and \eqref{Omr31}. The lemmas explain what should we do to get
equalities in \eqref{Omr2} and \eqref{Omr3}. In order to get
equality (at least up to terms $O(\de^3)$) in \eqref{Omr31}, we
have to leave only the first Fourier coefficient in
\eqref{Fourier2}. This leads to the following test-function:
\be\label{test-function}
\phi(x,y) = \left\{ \begin{array}{ll}
\cos\left(\frac{\pi a}{2}\right) e^{-\sqrt{M(a)}\, \delta (x-a)}
\sin\left(\frac{\pi y}{2}\right) & \textrm{$x \geq a$,}\\
\cos\left(\frac{\pi x}{2} \right) \sin\left(\frac{\pi y}{2}\right) &
\textrm{$x \leq a$.}
\end{array} \right.
\end{equation}

One can check that this function indeed gives us equality up to $O(\de^3)$
everywhere where it matters, and so it satisfies \eqref{Rq} which
proves \eqref{main1a}. The proof of the theorem in the case $a<1$
is thus finished.

\section{$a=1$}\label{a=1}
As before, the inequality \eqref{main2} is equivalent to the following two
inequalities:
\bee\label{maina1a}
\la\le \frac{\pi^2}{4}-M(a)\de^{2/3}+O(\de^{4/3})
\end{equation}
and
\bee\label{maina1b}
\la\ge \frac{\pi^2}{4}-M(a)\de^{2/3}-O(\de^{4/3}).
\end{equation}
We start by producing the test-function $\phi(x,y)$ with the
Rayleigh quotient equal to the RHS of \eqref{maina1a}. Such a
function is given by
\be\label{a1test}
\phi(x,y) = \left\{ \begin{array}{ll}
e^{-\pi^{2/3} \de^{1/3} (x-1)} \sin\left(\frac{\pi y}{2}\right) & \textrm{$x\geq 1$,}\\
\left(\cos(\mu )\right)^{-1} \cos(\mu  x )  \sin\left(\frac{\pi y}{2}\right) &
\textrm{$x\leq 1$,}
\end{array} \right.
\end{equation}
where $\mu = \frac{\pi}{2} - \pi^{1/3}\de^{2/3}$. A straightforward
(though rather lengthy) computation shows that
the Rayleigh quotient of this function indeed satisfies
\bee\label{a1quotient}
Q(\phi)=\frac{\pi^2}{4}-M(a)\de^{2/3}+O(\de^{4/3}),
\end{equation}
which proves \eqref{maina1a}. Another way of seeing that
\eqref{a1quotient} holds is to read the proof of
\eqref{maina1b} and check that the function $\phi$ changes all
inequalities in it into equalities.

Now we give a proof of \eqref{maina1b}. The proof is quite similar
to that of \eqref{main1b}. The biggest change is that instead of
Lemma \ref{left}, in the case $a=1$ we have the following result:

\begin{lem}\label{aequ1_L_Quadform1}
Let $f \in C^{2}[0,1]$ and  $f'(0)=0$. Then for all small enough positive $\de$
\be\label{aequ1_one_dim_aim}
\int_{0}^{1} \left(\frac{df}{dx}\right)^{2} dx - \left(\psf - M\de^{2/3} - C_{1}\de^{4/3}
\right) \int_{0}^{1} f^{2}(x) dx + \frac{\sqrt{M} f^{2}(1)}{2 \de^{2/3}} \geq 0,
\end{equation}
where
\be\label{aequ1_M}
M:= M(1) = \pi^{4/3}
\end{equation}
and $C_{1}$ is a constant, the precise value of which is not important.
Moreover, there exists another constant $C_{2}$, such that for
\be\label{aequ1_equality}
g(x) = \cos(\mu x)
\end{equation}
with $\mu^{2} = \psf - \pi^{4/3} \de^{2/3}$,
the inequality in the opposite direction is satisfied, namely:
\be\label{aequ1_one_dim_aim_inequality_in_other_direction}
\int_{0}^{1} \left(\frac{dg}{dx}\right)^{2} dx -
\left(\psf - M\de^{2/3} - C_{2}\de^{4/3} \right) \int_{0}^{1} g^{2}(x) dx +
\frac{\sqrt{M} g^{2}(1)}{2 \de^{2/3}} \leq 0.
\end{equation}
\end{lem}
\begin{proof}
The proof follows similar lines to those of Lemma \ref{left}.
Without loss of generality we can assume
that $f$ satisfies the extra boundary
condition
\bd
f'(1) = - \frac{\sqrt{M}\,f(1)}{2\de^{2/3}}.
\ed
Then
\bd
\int_{0}^{1} \left(\frac{df}{dx}\right)^{2} dx + \frac{\sqrt{M}\,f^{2}(1)}{2\de^{2/3}}
\ed
is the quadratic form of the operator $-\frac{d^{2}}{dx^{2}}$ with boundary conditions
\bd
f'(0)=0, \qquad f'(1)=-\frac{\sqrt{M}\,f(1)}{2 \de^{2/3}}.
\ed
We therefore need to prove that the first eigenvalue of
\be\label{aequ1_diff_equn}
- f''(x)  = \mu^{2} f(x), \qquad f'(0)=0, \qquad f'(1)=-\frac{\sqrt{M}\, f(1)}{2\de^{2/3}}
\end{equation}
equals $\psf- M \de^{2/3} + O(\de^{4/3})$ (this indeed would prove both
statements of lemma).  The eigenvalues of \eqref{aequ1_diff_equn}
are the values of $\mu_j^{2}$ where $\mu_j$ are solutions to
\be\label{aequ1_sqrtM_tan}
\mu \tan(\mu) = \frac{\sqrt{M}}{2 \de^{2/3}};
\end{equation}
the corresponding eigenfunctions are $\cos(\mu_j x)$. Therefore, we need to
make sure that
\be\label{M11}
\sqrt{\psf- M \de^{2/3}}\tan\left(\sqrt{\psf- M \de^{2/3}}\right)=\frac{\sqrt{M}}{2
\de^{2/3}}+O(1).
\end{equation}
The LHS of \eqref{M11} is
\begin{equation*}
\frac{\pi}{2}\frac{\pi}{M\de^{2/3}}+O(1)=\frac{\pi^2}{2M\de^{2/3}}+O(1).
\end{equation*}
This is precisely the RHS of \eqref{M11} iff $M$ is given by
\eqref{aequ1_M}.
\end{proof}

We will now prove \eqref{maina1b}.  To do this it is enough to show that for arbitrary
$u \in C^{\infty}(\Om')$ such that
$u(x,0)=\frac{\partial u}{\partial y}(x,1) =\frac{\partial u}
{\partial x}(0,y)= 0$, we have
\be\label{a1quadratic}
\int\!\!\!\int_{\Om'} \vert \nabla u \vert^{2} dxdy - \left( \psf -M \de^{2/3} - O(\de^{4/3})
\right) \int\!\!\!\int_{\Om'} \vert u \vert^{2} dxdy \geq 0.
\end{equation}
The LHS of \eqref{a1quadratic} can be
rewritten as
\bd
\bs
& \left(\int\!\!\!\int_{\Omega_{l}} \vert \nabla u \vert^{2} dxdy -
\left( \psf -M \de^{2/3} - O(\de^{4/3}) \right) \int\!\!\!\int_{\Omega_{l}}
\vert u \vert^{2} dxdy\right)\\
+& \left(\int\!\!\!\int_{\Omega_{r}} \vert \nabla u \vert^{2} dxdy -
\left( \psf -M \de^{2/3} - O(\de^{4/3}) \right) \int\!\!\!\int_{\Omega_{r}}
\vert u \vert^{2} dxdy\right)\\
=&: L(u) + R(u).
\end{split}
\ed
Similar to Section \ref{a<1}, we will show that
\be\label{aequ1_bound_R}
R(u) \geq \frac{\sqrt{M}}{2\de^{2/3}} \int_{1-\de}^{1} \vert u(1,y)\vert^{2} dy
\end{equation}
and
\be\label{aequ1_bound_L}
L(u) +  \frac{\sqrt{M}}{2\de^{2/3}} \int_{1-\de}^{1} \vert u(1,y)\vert^{2} dy \geq 0.
\end{equation}
The proof of  \eqref{aequ1_bound_R} is absolutely analogous to the
proof of \eqref{Omr1}, and we will skip it. In order to prove
\eqref{aequ1_bound_L}, it is sufficient to show that
\be\label{lemmay}
\begin{split}
\int\!\!\!\int_{\Omega_{l}} \left(\frac{\partial u}{\partial x} \right)^{2} dxdy &-
\left( \frac{\pi^{2}}{4} - M \de^{2/3} - C_{1}\de^{4/3} \right)
\int\!\!\!\int_{\Omega_{l}} |u|^{2} dxdy\\ &+ \frac{\sqrt{M}}{2\de^{2/3}} \int_{1-\delta}^{1}
\vert u (1,y) \vert^{2} dy \geq 0
\end{split}
\end{equation}
for some choice of the constant $C_{1}$.  This follows immediately if we
apply Lemma \ref{aequ1_L_Quadform1} to the function $u(\cdot,y)$ and then integrate
over $y$.  The proof of \eqref{maina1b} is therefore complete. The
careful look at the proof together with the second part of
Lemma \ref{aequ1_L_Quadform1}  shows that the choice of the
test-function \eqref{a1test} indeed changes all the inequalities
into equalities, and so \eqref{maina1a} is proved. This finishes
the proof of our theorem in the case $a=1$.

\section{Arbitrary $a$}
Consider now the case of arbitrary $a$. The main
difference between this case and the case $a\le 1$ is the fact that
now we have to take care of several test-functions, using the
mini-max principle. The cases of integer and non-integer $a$
require slightly different approach as well as the cases of even
and odd $[a]$. We consider in details the case of non-integer $a$
with even integer part and prove the theorem in this case. Proof of the
other cases is similar, and we will not give it
here. So, let us assume that $[a]=2l$, $l\in\N$,
$a\not\in\N$.
Then, as we have mentioned already, the number of eigenvalues
$n=2l+1$, and the top
eigenvalue $\la_{2l+1}$ comes from the Neumann problem
$L_N$. This problem has exactly
$l+1$ eigenvalues (the other $l$ eigenvalues are due to $L_D$).
As before, we construct the test-functions to prove the upper
bound
\bee\label{main2a}
\la_{2l+1}\le \frac{\pi^2}{4}-M(a)\de^2+O(\de^3)
\end{equation}
and estimate the quadratic form to prove the lower bound
\bee\label{main2b}
\la_{2l+1}\ge \frac{\pi^2}{4}-M(a)\de^2-O(\de^3).
\end{equation}
More precisely, in order to prove \eqref{main2a} we will construct
$l+1$ functions $\phi_1,\dots,\phi_{l+1}\in H^1(\Om')$,
$\phi_j(x,0)=\frac{\di \phi_j}{\di
x}(0,y)=\frac{\di \phi_j}{\di y}(x,1)=0$
such that every non-trivial linear combination
$\phi=\sum_{j=1}^{l+1}\al_j\phi_j$ satisfies \eqref{Rq}. As in
section \ref{a<1}, it is relatively easy to construct
test-functions which satisfy \eqref{Rq} with a weaker constant.
So, once again we will start with proving the lower bound and this
proof will give us the recipe for choosing the optimal test-functions.
Taking into account that we are estimating the eigenvalue number
$l+1$, the variational formulation of \eqref{main2b} is the
following: for every set
$u_1,\dots,u_{l+1}\in C^\infty(\Om')$,
$u_j(x,0)=0$, $\frac{\di u_j}{\di
x}(0,y)=0$, $\frac{\di u_j}{\di y}(x,1)=0$ there exist constants
$\al_j$, not all zeros ($1\le j\le l+1$), such that if
$u=\sum_{j=1}^{l+1}\al_j u_j$, then
$Q(u)\ge \frac{\pi^2}{4}-M(a)\de^2-O(\de^4)$. The proof of this
statement goes similarly to the proof from section \ref{a<1}, but
we need to make two modifications. First of all,
instead of lemma \ref{left} we use the following lemma:

\begin{lem}\label{nonint_L_N_Quadform1}
Let $2l<a<2l+1, l \in \N$ and let  $f_{1},\dots, f_{l+1} \in C^{2}[0,a]$ satisfy
$f_{j}'(0)=0$. Then there exist constants $\alpha_{j}$ not all zeros ($1 \leq j \leq l+1$),
such that for the linear combination $f=\sum_{j=1}^{l+1} \alpha_{j} f_{j}$ the following
inequality is satisfied:
\be\label{nonint_N_one_dim_aim}
\int_{0}^{a} \left( \frac{df}{dx}\right)^{2}dx - \psf \int_{0}^{a} f^{2}(x) dx +
\frac{\sqrt{M(a)}\, f^{2}(a)}{2} \geq 0,
\end{equation}
where
\be\label{nonint_N_M(a)}
M(a) = \ps \tan^{2} \left( \frac{\pi a}{2}\right).
\end{equation}
Moreover, inequality in the other direction is reached in \eqref{nonint_N_one_dim_aim}
for any linear combination of the following functions:
\be\label{nonint_N_equality}
f_{j}(x) = \cos(\mu_{j} x ), \qquad 1\leq j \leq l+1,
\end{equation}
where $\mu_{j}$ are the first $(l+1)$ solutions to
 \be\label{nonint_N_laj_equn}
\mu_{j} \tan(\mu_{j} a ) = \frac{\pi}{2} \tan\left(\frac{\pi a}{2}\right),
\end{equation}
with $\mu_{l+1} = \frac{\pi}{2}$.
\end{lem}
\begin{proof}
The proof of this lemma is similar to the proof of Lemma
\ref{left}. First of all we notice that without loss of generality
we can assume that $f_j$ satisfy the additional boundary condition
$\qquad f_j'(a)=-\frac{\sqrt{M(a)}\,f_j(a)}{2}$. Then the
statement of the lemma is equivalent to the fact that the
$(l+1)$st eigenvalue of the problem
\be\label{nonint_N_diff_equn}
-f''(x)=\mu^{2} f(x), \qquad f'(0)=0, \qquad f'(a)=-\frac{\sqrt{M(a)}\,f(a)}{2}
\end{equation}
is $\psf$ and the functions \eqref{nonint_N_equality} are the
first $(l+1)$ eigenfunctions. As in the proof of Lemma \ref{left}
the eigenvalues of \eqref{nonint_N_diff_equn} are the values of $\mu_j^{2}$,
where $\mu_j$ are (positive) solutions to
\be\label{nonint_M(a)_tan}
\mu_j \tan(\mu_j a) = \frac{\sqrt{M(a)}}{2}=
\frac{\pi}{2} \tan\left(\frac{\pi a}{2}\right), \qquad
\mu_1<\mu_2<\dots
\end{equation}
Obviously, $\mu_j=\pi/2$ solves \eqref{nonint_M(a)_tan}, and we
just have to find the number of solutions of
\eqref{nonint_M(a)_tan} which are smaller than $\ps$. It is easy to see
that \eqref{nonint_M(a)_tan} has precisely one solution in the interval
$(0,\frac{\pi}{2a})$ and in each of the intervals
$(\frac{(2j-1)\pi}{2a},\frac{(2j+1)\pi}{2a})$ ($j=1,2,\dots$).
Since $\pi/2\in (\frac{(2l-1)\pi}{2a},\frac{(2l+1)\pi}{2a})$, we
see that indeed $\mu_{l+1}=\pi/2$. This finishes the proof of
lemma.
\end{proof}

There is another problem arising if one tries to use the proof
from the previous sections directly to prove the inequality
\eqref{Omr2} (the proof of \eqref{Omr1} is unchanged)
We can not apply lemma \ref{nonint_L_N_Quadform1}
to the set of
functions $u_j(\cdot,y)$ for each $y$ separately and then integrate
over $y$ (like we did
in the previous sections), since the choice
of coefficients $\al_j$ would then depend on $y$. Therefore, we have to
use the Fourier decomposition in the region above the obstacle to
show that, roughly speaking, only the first Fourier term matters,
which would allow us to apply lemma \ref{nonint_L_N_Quadform1}
only to this first term. So, let us write $u$ in
terms of a Fourier Series for $0\le x\leq a$:
\bd
u(x,y) = \sum_{j=0}^{\infty} \hat{u}_{j}(x) \hat{\eta}_{j}(y),
\ed
where
\bd
\hat{\eta}_{j}(y) = \left\{ \begin{array}{ll}
\de^{-1/2} & \textrm{$j=0$,}\\
\sqrt{\frac{2}{\de}} \, \cos\left(\frac{\pi j}{\de} (1-y)\right) & \textrm{$ j \geq 1$.}
\end{array} \right.
\ed
After simple calculations we have
\be\label{nonint_nabla_u_FS}
\bs
\int\!\!\!\int_{\Om_{l}} \vert \nabla u \vert^{2} dxdy
&=  \int_{0}^{a} \left( \hat{u}_{0}'(x) \right)^{2} dx + \sum_{j=1}^{\infty}\int_{0}^{a}
\left( \hat{u}_{j}'(x) \right)^{2} dx\\ &+ \frac{\ps}{\de^{2}} \sum_{j=1}^{\infty} j^{2}
\int_{0}^{a} \hat{u}_{j}^{2}(x) dx\\
\geq \int_{0}^{a} \left( \hat{u}_{0}'(x) \right)^{2} dx &+ \sum_{j=1}^{\infty}\int_{0}^{a}
\left( \hat{u}_{j}'(x) \right)^{2} dx + \frac{\ps}{\de^{2}} \sum_{j=1}^{\infty}
\int_{0}^{a} \hat{u}_{j}^{2}(x) dx
\end{split}
\end{equation}
and
\be\label{nonint_u_FS}
\int\!\!\!\int_{\Om_{l}} \vert u \vert^{2} dxdy = \int_{0}^{a} \hat{u}_{0}^{2}(x) dx +
\sum_{j=1}^{\infty} \int_{0}^{a} \hat{u}_{j}^{2}(x) dx.
\end{equation}
Recall that we are trying to prove \eqref{Omr2}, i.e. that
\be\label{nonint_apply_lem_to_this}
\int\!\!\!\int_{\Om_{l}} \vert \nabla u  \vert^{2} dxdy - \psf \int\!\!\!\int_{\Om_{l}}
\vert u \vert^{2} dxdy + \frac{\sqrt{M(a)}}{2} \int_{1-\de}^{1} \vert u(a,y)\vert^{2} dy
\geq 0.
\end{equation}
Substituting \eqref{nonint_nabla_u_FS} and \eqref{nonint_u_FS} into
the LHS of \eqref{nonint_apply_lem_to_this} gives
\be\label{Fourier17}
\bs
&\int\!\!\!\int_{\Om_{l}} \vert \nabla u  \vert^{2} dxdy - \psf \int\!\!\!\int_{\Om_{l}}
\vert u \vert^{2} dxdy + \frac{\sqrt{M(a)}}{2} \int_{1-\de}^{1} \vert u(a,y)\vert^{2} dy\\
&\geq \left( \int_{0}^{a} \left( \hat{u}_{0}'(x)\right)^{2} dx - \psf \int_{0}^{a}
\hat{u}_{0}^{2}(x) dx + \frac{\sqrt{M(a)}\, \hat{u}_{0}^{2}(a)}{2} \right)\\
& + \sum_{j=1}^{\infty} \left( \int_{0}^{a} \left( \hat{u}_{j}'(x)\right)^{2} dx + \ps
\left(\frac{1}{\de^{2}}-\frac{1}{4} \right) \int_{0}^{a}  \hat{u}_{j}^{2}(x) dx +
\frac{\sqrt{M(a)}\, \hat{u}_{j}^{2}(a)}{2}\right) \\
&\geq  \int_{0}^{a} \left( \hat{u}_{0}'(x)\right)^{2} dx - \psf \int_{0}^{a}
\hat{u}_{0}^{2}(x) dx + \frac{\sqrt{M(a)}\, \hat{u}_{0}^{2}(a)}{2}.
\end{split}
\end{equation}
The aim of this exercise was to reduce proving the
two-dimensional inequality \eqref{nonint_apply_lem_to_this} to
the one-dimensional one
\bee\label{onedim}
 \int_{0}^{a} \left( \hat{u}_{0}'(x)\right)^{2} dx - \psf \int_{0}^{a}
\hat{u}_{0}^{2}(x) dx + \frac{\sqrt{M(a)}\,
\hat{u}_{0}^{2}(a)}{2}\ge 0,
\end{equation}
and this is precisely what Lemma \ref{nonint_L_N_Quadform1} is
about. Namely, having a set of $(l+1)$ functions
$u_1$,...,$u_{l+1}$, we apply Lemma \ref{nonint_L_N_Quadform1} to
their first Fourier coefficients to get a linear combination of
them which satisfies \eqref{onedim}. Computations above show that
this implies \eqref{Omr2} and \eqref{main2b}.

In order to prove \eqref{main2a}, we, as always, just look
carefully at the proof and try to find functions which
change all the inequalities into
equalities (at least up to $O(\de^3)$ terms). This leads to the
following set of test-functions:
\be\label{nonint_N_optimal_test_functions}
\phi_{j}(x,y) = \left\{ \begin{array}{ll}
f_{j}(a) e^{-\sqrt{M(a)}\, \delta (x-a)} \sin\left(\frac{\pi y}{2}\right) &
\textrm{$x \geq a$,}\\
f_{j}(x) \sin\left(\frac{\pi y}{2}\right) & \textrm{$x \leq a$,}
\end{array} \right.
\end{equation}
where $M(a)$ is as in \eqref{nonint_N_M(a)} and $f_{j}(x)$ are
\eqref{nonint_N_equality}. It is an easy matter to check that any linear combination
of them $\phi=\sum\al_j\phi_j$ satisfies
\bee\label{Rq2}
Q(\phi):=\frac{\int_{\Om'}|\nabla\phi|^2 dxdy}{\int_{\Om'}|\phi|^2 dxdy}\le
\frac{\pi^2}{4}-M(a)\de^2+O(\de^3),
\end{equation}
which proves \eqref{main2a}. This finishes the proof of the
theorem in the case $a\ne [a]=2l$. Other cases are treated
similarly.

\end{document}